\input amstex
\documentstyle{amsppt}
\magnification=\magstep1 \hsize=5.2in \vsize=6.7in

\catcode`\@=11

\topmatter

\title On $\text{\rm II}_1$ factors arising from 2-cocycles of $w$-rigid
groups
\endtitle
\author{Remus Nicoara$^1$, Sorin Popa, Roman Sasyk}
\endauthor
\leftheadtext{Remus Nicoara, Sorin Popa, Roman Sasyk}
\thanks
$^1$ supported by NSF under Grant No. DMS 0500933
\endthanks

\abstract We consider $\text{\rm II}_1$ factors $L_\mu(G)$ arising
from $2$-cocyles $\mu \in \text{\rm H}^2(G,\Bbb T)$ on groups $G$
containing infinite normal subgroups $H \subset G$ with the relative
property $\text{\rm(T)}$ (i.e. $G$ {\it w-rigid}).
 We prove that given any separable $\text{\rm II}_1$ factor $M$,
the set of $2$-cocycles $\mu_{|H}\in \text{\rm H}^2(H,\Bbb T)$ with
the property that $L_\mu(G)$ is embeddable into $M$ is at most
countable. We use this result, the relative property (T) of $\Bbb
Z^2 \subset \Bbb Z^2 \rtimes \Gamma$ for $\Gamma \subset SL(2,\Bbb
Z)$ non-amenable and the fact that every cocycle $\mu_\alpha \in
{\text{\rm H}}^2(\Bbb Z^2,\Bbb T)\simeq \Bbb T$ extends to a cocycle
on $\Bbb Z^2 \rtimes SL(2,\Bbb Z)$, to show that the one parameter
family of II$_1$ factors $M_\alpha(\Gamma)=L_{\mu_{\alpha}}(\Bbb Z^2
\rtimes \Gamma)$, $\alpha \in \Bbb T$, are mutually non-isomorphic,
modulo countable sets, and cannot all be embedded into the same
separable II$_1$ factor. Other examples and applications are
discussed.
\endabstract

\endtopmatter

\document

\heading {0. Introduction} \endheading

Ever since Connes' celebrated ``rigidity'' paper ([C]), groups with
the property $\text{\rm(T)}$ of Kazhdan have played an important
r\^{o}le in operator algebra, being used to obtain a plethora of
rigidity results and interesting examples (see e.g. [Ch1,2], [O1],
[P1-4], [V]), especially in the theory of II$_1$ factors. More
recently, a weaker version of the property $\text{\rm(T)}$, merely
requiring the existence of a ``large'' subgroup with the relative
property (T) of Kazhdan-Margulis ([M], [K]), proved to be equally
important (cf. [P2,3,4]). 
The prototype of such group is $\Bbb Z^2 \rtimes
SL(2,\Bbb Z)$, with $\Bbb Z^2$ its relative property $\text{\rm(T)}$
subgroup (cf. [K], [M]). Thus, it is shown in ([P2]) that the
$\text{\rm II}_1$ factors associated with this arithmetic group, and
more generally with the groups $\Bbb Z^2\rtimes\Gamma$, for $\Gamma$
non amenable finitely generated subgroups of $SL(2,\Bbb Z)$, have
trivial fundamental group and are non isomorphic if the groups
$\Gamma$ have different $\ell^2$-Betti numbers, e.g. $\Gamma=\Bbb
F_n, \,n=2,3,\dots$. (For $\Gamma \subset SL(2,\Bbb Z)$ non-amenable
the inclusion of groups $\Bbb Z^2\subset\Bbb Z^2\rtimes\Gamma$ was
shown to have relative property (T) in \cite{Bu}). This provided the
first examples of factors with trivial fundamental group ([MvN]).

More generally, in ([P2]; see also [P4]) one considers a one
parameter family of II$_1$ factors $M_\alpha(\Gamma), \alpha \in
\Bbb T,$ associated with $\Bbb Z^2 \rtimes \Gamma$, for each $\Gamma
\subset SL(2,\Bbb Z)$ non-amenable, and one proves several rigidity
properties and classification results for $M_\alpha(\Gamma)$. We
continue in this paper the analysis of this interesting class of
II$_1$ factors.

The factors $M_\alpha(\Gamma)$ are defined to be crossed product
$\text{\rm II}_1$ factors of the form
$M_{\alpha}(\Gamma)=R_{\alpha}\rtimes_{\sigma_{\alpha}} \Gamma$,
where $\alpha \in \Bbb T$, $R_\alpha$ is the finite von Neumann
algebra generated by two unitaries $u,v \in R_\alpha$ satisfying the
relation $uv=\alpha vu$ and trace $\tau(u^kv^l)=0, \forall (k,l)\neq
(0,0)$, $\Gamma$ is an arbitrary non-amenable subgroup of
$SL(2,\Bbb{Z})$ and the action $\sigma_{\alpha}$ is implemented by
the restriction to $\Gamma$ of the action of $SL(2,\Bbb{Z})$ on
$R_\alpha$ given by
$\sigma_{\alpha}(g)(u^kv^l)=
\alpha^{\frac{1}{2}(kl-(ak+bl)(ck+dl))}u^{ak+bl}v^{ck+dl}$
 where $g=\pmatrix a & b \cr c & d \cr
\endpmatrix \in \Gamma$ (see [B]).

If $\alpha$ is a primitive root of unity of order $n$, then
$R_\alpha$ is isomorphic to $L((n\Bbb Z)^2)\otimes M_{n\times
n}(\Bbb C)$  and $M_\alpha(\Gamma)\simeq L((n\Bbb
Z)^2\rtimes\Gamma)\otimes M_{n\times n}(\Bbb C)$ (Corollary 5.2.1 of
[P2]). If $\Gamma$ is finitely generated and $\alpha'$ is another
primitive root of 1 of order $n'$ then by ([P2]) $M_{\alpha}(\Gamma)
\simeq M_{\alpha'}(\Gamma)$ if and only if $n=n'$. If in turn
$\alpha=e^{2\pi i\theta} \in\Bbb T$ with $\theta\in [0,1/2)$
irrational then $R_\alpha$ is isomorphic to the hyperfinite II$_1$
factor, represented as the {\it irrational rotation} von Neumann
algebra $R_\alpha$ ([R]). The factors $M_\alpha(\Gamma)$ are called
{\it irrational} (resp. {\it rational}) {\it rotation} HT {\it
factors} when $\alpha = e^{2\pi i\theta}$ with $\theta \in [0,1/2)
\setminus \Bbb Q$ (resp $\theta \in \Bbb Q$). By ([P2]), if $\Gamma$
is non-amenable then an irrational rotation HT factor
$M_\alpha(\Gamma)$ cannot be embedded into a rational rotation HT
factor $M_{\alpha'}(\Gamma')$.

The problem of classifying the family of factors $M_\alpha(\Gamma)$,
in terms of the embedding $\Gamma \subset SL(2,\Bbb Z)$ and the
parameter $\alpha \in \Bbb T$, is quite natural. In this respect, it
has been conjectured in ([P2]) that for each fixed $\Gamma \subset
SL(2,\Bbb Z)$ non-amenable (notably for $\Gamma=SL(2,\Bbb Z)$), the
factors $M_\alpha(\Gamma)$, $\alpha \in \Bbb T$ irrational, are
mutually non-isomorphic. In this paper we will give a partial,
positive answer to this problem, by showing that for each
fixed non-amenable group $\Gamma \subset SL(2,\Bbb Z)$ the factors
$M_\alpha(\Gamma)$, $\alpha \in \Bbb T$, are mutually non stably
isomorphic, modulo countable sets, i.e. there 
are at most countably many
$\alpha$'s in $\Bbb T$ such that
$M_\alpha(\Gamma)\simeq M_{\alpha_0}(\Gamma)$,
for a fixed, arbitrary $\alpha_0\in \Bbb T$.

We will alternatively view a factor $M_\alpha(\Gamma)$ as a cocycle
group von Neumann algebra $L_{\mu_\alpha} (\Bbb Z^2 \rtimes \Gamma)$
(see [CJ]) corresponding to a projective left regular representation
$\lambda_{\mu_\alpha}$ with the scalar 2-cocycle $\mu_\alpha \in
\text{\rm H}^2(\Bbb Z^2 \rtimes \Gamma, \Bbb T)$ depending on
$\alpha\in\Bbb T$. To explain this, let us first recall some
definitions.

Let $G$ be a discrete group and $\mu \in \text{\rm H}^2(G,\Bbb T)$ a
2-cocycle on $G$, i.e. $\mu: G \times G \rightarrow \Bbb T$
satisfies $\mu_{g,h}\mu_{gh,k} = \mu_{h,k}\mu_{g,hk}$, $\forall
g,h,k\in G$. One associates to $\mu$ the {\it projective left
regular representation} $\lambda_\mu : G \rightarrow \Cal
U(l^2(G))$, defined by $\lambda_\mu(g) (\Sigma_{h\in G} c_h
\xi_h)=\Sigma_{h\in G} c_h \mu_{g,h} \xi_{gh}$, where
$\{\xi_h\}_{h\in G}$ is the canonical basis of $l^2(G)$. Denote by
$L_\mu(G)=\lambda_\mu(G)''$ the cocycle group von Neumann algebra of
$(G,\mu)$. It is well known that one has an isomorphism $\text{\rm
H}^2(\Bbb Z^2,\Bbb T)\cong\Bbb T$, taking $\alpha \in \Bbb T$ to
$\mu_\alpha\in\text{\rm H}^2(\Bbb Z^2,\Bbb T)$, where:
$$\mu_\alpha((k,l),(k',l'))=\alpha^{\frac{1}{2}({kl'-k'l})}$$

If we define $R_{\mu_{\alpha}}$ to be the cocycle group von Neumann
algebra $L_{\mu_\alpha}(\Bbb{Z}^2)$, then $R_{\mu_{\alpha}}$ is
generated by the unitary elements $u=\lambda_{\mu_{\alpha}}(1,0)$,
$v=\lambda_{\mu_{\alpha}}(0,1)$, which satisfy the relation
$uv=\alpha vu$, thus being naturally isomorphic to $R_\alpha$.
Moreover, $\mu_\alpha$ is invariant to the action $\sigma$ of
$SL(2,\Bbb{Z})$ on $\Bbb{Z}^2$, thus $\sigma$ implements an action
$\sigma_{\mu_{\alpha}}$ of $SL(2,\Bbb{Z})$ on
$R_{\mu_{\alpha}}=R_\alpha$ which coincides with the action
$\sigma_\alpha$ defined above.

Since any $\Gamma\subset SL(2,\Bbb Z)$ has Haagerup's compact
approximation property ([H]), by (6.9.1 in [P2]) it follows that
$M_\alpha(\Gamma)$ has Haagerup's property relative to $R_\alpha$
(as defined in 2.1 of [P2]). Also, by (Example 2 on page 62 of [Bu])
the pair $(\Bbb{Z}^2 \rtimes \Gamma, \Bbb{Z}^2)$ has the relative
property (T) for any non-amenable subgroup $\Gamma \subset SL(2,
\Bbb Z)$ and thus, by (6.9.1 in [P2]), the embedding $R_\alpha
\subset M_\alpha(\Gamma)$ is rigid in the sense of (Definition 4.2
in [P2]).

Since the action of $SL(2, \Bbb Z)$ on $\Bbb Z^2$ is outer, by
(3.3.2 $(ii)$ in [P2]) $\sigma_\alpha$ are properly outer actions of
$SL(2, \Bbb Z)$ (thus of $\Gamma$ as well) on $R_\alpha$.
Furthermore, since the stabilizer of any non-trivial element in
$\Bbb Z^2$ is a cyclic group, it follows that if $\Gamma$ leaves a
finite subset $\neq \{(0,0)\}$ of $\Bbb Z^2$ invariant, then it is
almost cyclic. Thus, by (3.3.2 $(i)$ in [P2]) any non-amenable
$\Gamma \subset SL(2, \Bbb Z)$ acts ergodically on $R_\alpha$. Thus,
$R_\alpha \subset M_\alpha(\Gamma)$ satisfies $R_\alpha'\cap
M_\alpha(\Gamma)\subset R_\alpha$. In particular, when $\alpha$ is
irrational, $R_\alpha \subset M_\alpha(\Gamma)$ are irreducible
inclusions of $\text{\rm II}_1$ factors, and they are HT
inclusions in the
sense of (6.1 in [P2]).

By ([P2]) the factors $M_\alpha(\Gamma)$ are non-$\Gamma$ and by
([O2]) they are prime, i.e., they cannot be decomposed into a tensor
product of $\text{\rm II}_1$ factors. It was shown in [P4] that two
factors $M_\alpha(\Gamma)$ with $\Gamma$ torsion free are isomorphic
iff $\sigma_\alpha(\Gamma)$ are cocycle conjugate in Out$(R)$. In
particular, isomorphism between irrational rotation HT factors
$M_\alpha(\Gamma)$, with torsion free $\Gamma$, implies isomorphism
of the corresponding groups $\Gamma$. Also, it follows from ([P2])
that if $\Gamma$ is torsion free then $M_\alpha(\Gamma)$ has
countable fundamental group (see the Appendix for a more general
result).

The factors $M_\alpha(\Gamma)$ are easily seen to be ``approximately
embeddable'' into the hyperfinite $\text{\rm II}_1$ factor $R$ (in
the sense of Connes [C1]), i.e., $M_\alpha(\Gamma) \subset
R^\omega$. Indeed, let $m_k/n_k$ be a sequence of rational numbers
such that $\alpha_k=exp(2\pi i m_k/n_k) \rightarrow \alpha$ and $n_k
\rightarrow \infty$. Let $\pi_k$ be the projective representation
with 2-cocycle $\mu_{\alpha_k}$, of the group $\Bbb Z^2 \rtimes
SL(2, \Bbb Z)$ on $\Cal H_n=\ell^2((\Bbb Z/n_k\Bbb Z)^2\rtimes SL(2,
\Bbb Z/n_k\Bbb Z))$. Then $g \mapsto (\pi_n(g))_n$ is an embedding
of $L_\alpha(\Bbb Z^2) \rtimes SL(2, \Bbb Z)$ into $\Pi_n \Cal
B(\Cal H_n)\subset R^\omega$. However, we have:

\proclaim {0.1. Theorem} Let $M$ be a separable $\text{\rm II}_1$
factor. For each fixed $\Gamma\subset SL(2,\Bbb Z)$ non-amenable
there exist at most countably many
$\alpha \in \Bbb T$ such that $M_{\alpha}(\Gamma)$ is embeddable
into $M$ $($not necessarily unitaly$)$. In particular, the factors
$\{ M_{\alpha}(\Gamma)\} _{\alpha \in \Bbb T}$ are non stably
isomorphic modulo countable sets.
\endproclaim

Note that the above theorem gives an alternative proof to Ozawa's
result on the non-existence of universal separable $\text{\rm II}_1$
factors in [O1], without using Gromov's property $\text{\rm(T)}$
groups. More precisely, in the same spirit as the results in [O1],
the above theorem shows that there exist no separable finite von
Neumann algebra $M$ that can contain an uncountable set of
projective unitary representations $\{\pi_j\}_{j\in J}$ of
$\Bbb{Z}^2\rtimes SL(2, \Bbb{Z})$ with distinct scalar 2-cocycles
$\{\mu_{\pi_j}\}_j$. Theorem 0.1 will follow  as a special case of
the following:

\proclaim{0.2. Theorem} Let $H\subset G$ be an inclusion of
discrete groups with the relative property $\text{\rm (T)}$. Let $M$
be a separable finite von Neumann algebra. Let $J$ be the set of
scalar $2$-cocycles $\mu \in \text{\rm H}^2(G,\Bbb T)$ such that
$L_\mu(G)$ can be embedded into $M$ $($not necessarily unitaly$)$.
Then the set $\{\mu_{|H} \mid \mu \in J\}\subset \text{\rm
H}^2(H,\Bbb T)$ is countable.
\endproclaim

We prove this result in Section 1, by using a separability argument
similar to ([C], [P1], [P2], [O1], [GP]) and a characterization of
the relative property $\text{\rm(T)}$ in terms of projective
representations. In Section 2 we give examples of pairs of groups
$H\subset G$  with the relative property $\text{\rm (T)}$ with the
torus $\Bbb T$ embedded as a subgroup of $2$-cocycles $\Bbb T
\subset \text{\rm H}^2(G,\Bbb T)$, such that $\Bbb T \ni \mu \mapsto
\mu_{|H} \in \text{\rm H}^2(H,\Bbb T)$ is one to one. In Section 3
we give an explicit description of the disintegration of type II$_1$
von Neumann algebras from the property (T) groups $\Lambda$
constructed by Serre (see \cite{dHVa} pg. 40), as central extensions
of property (T) groups $\Gamma$, in terms of factors
$L_\alpha(\Gamma)$ associated with 2-cocycles of $\Gamma$. We also
show that the factors in the
disintegration of the algebra $L(\Lambda)$ of an arbitrary property
(T) group $\Lambda$ are mutually non-isomorphic, modulo countable
sets, by using a separability argument similar to Ozawa's proofs in ([O1]).

There are strong indications from ([P2]) and results in this paper
that for each $\Gamma\subset SL(2,\Bbb Z)$ non-amenable the factors
$\{ M_{\alpha}(\Gamma)\} _{\alpha\in I}$, where $I=\{e^{2\pi it}\mid
t \in [0,1/2)\setminus \Bbb Q\}$, are all mutually non stably
isomorphic and have trivial fundamental group, and that if the
normalizer of $\Gamma$ in $GL(2, \Bbb Z)$ is equal to $\Gamma$ then
$\text{\rm Out}(M_\alpha(\Gamma))$ is isomorphic to the character
group of $\Gamma$. Furthermore, if $\Gamma_1,\Gamma_2 \subset
SL(2,\Bbb Z)$ non-amenable and $\alpha_1, \alpha_2 \in I$  then
$M_{\alpha_1}(\Gamma_1)$ should be isomorphic to
$M_{\alpha_2}(\Gamma_2)$ if and only if $\alpha_1 = \alpha_2$ and
$\Gamma_1, \Gamma_2$ are conjugate in $SL(2, \Bbb Z)$, by an
automorphism of $SL(2, \Bbb Z)$ (thus, by an element in $GL(2, \Bbb
Z)$). In this respect, note that by ([P4]) the non-isomorphism of
the factors $M_{\alpha}(\Gamma)$ for a fixed $\Gamma$ amounts to
showing that if $\alpha_1\neq\alpha_2$ then
$\sigma_{\alpha_1},\sigma_{\alpha_2}$ are not cocycle conjugate.
While we cannot prove this fact, we obtain in Section 4 the
following result, whose proof is inspired from an argument in [Ch]:

\proclaim {0.3. Theorem} Let $\Gamma \subset SL(2,\Bbb{Z})$ be a
subgroup of $SL(2,\Bbb{Z})$ containing a parabolic element $a$ and
an element $b$ that doesn't commute with $a$. If for some
$\alpha_1,\alpha_2 \in I$ the actions $\sigma_{\alpha_i}$ of
$\Gamma$ on $R_{\alpha_i}, i=1,2,$ are conjugate, then
$\alpha_1=\alpha_2$.
\endproclaim

An initial version of this paper has been circulated since January
2004 under the title ``Some remarks on irrational rotation HT
factors'' math.OA/0401139. It is a pleasure for us to thank Bachir
Bekka and Larry Brown for comments and useful discussions on the
first version of the paper.

Part of this work was done while R.N. and R.S. were visiting UCLA.
They would like to thank UCLA for the hospitality.

\heading 1. A cocycle characterization of relative property $\text{\rm(T)}$
\endheading

Recall that an inclusion of discrete groups $H \subset G$ has the
{\it relative property} $\text{\rm(T)}$ of Kazhdan-Margulis if there
exist $\delta_0>0$ and a finite subset $F_0\subset G$ such that if
$\pi:G\rightarrow\Cal{U(H)}$ is a unitary representation of $G$ on
the Hilbert space $\Cal H$ with a unit vector $\xi\in\Cal H$
satisfying $\|\pi(g)\xi-\xi\|<\delta_0$ for all $g\in F_0$, then
there exists a vector $\xi_0\in \Cal H$ such that $\pi(h)\xi_0=\xi_0$
for all $h\in H$. Note that in case $H=G$ this amounts to $G$ itself
having Kazhdan's property $\text{\rm(T)}$. It is easy to see that if
$H$ is normal in $G$ then the above condition is equivalent to the
following:

\vskip .05in \noindent {\bf (1.0).} $\forall \varepsilon>0$ there
exist a finite subset $F(\varepsilon)\subset G$ and
$\delta(\varepsilon)>0$ such that if $\pi:G\rightarrow\Cal {U(H)}$
is a unitary representation of $G$ on the Hilbert space $\Cal H$
with a unit vector $\xi\in\Cal H$ satisfying
$\|\pi(g)\xi-\xi\|<\delta(\varepsilon)$, $\forall g\in
F(\varepsilon)$, then there exists a unit vector $\xi_0\in \Cal H$
such that $\|\xi_0-\xi\|<\varepsilon$ and $\pi(h)\xi_0=\xi_0$,
$\forall h\in H$. \vskip .05in

The above condition is in fact equivalent to the relative property
$\text{\rm(T)}$ even for inclusions that are not necessarily normal
(cf. \cite{J}), but we will only use here the equivalence for normal
inclusions.

We first show that if $H\subset G$ has the relative property
$\text{\rm(T)}$ then the projective representations satisfy a
property similar to $(1.0)$. To state the result, recall that a
projective (unitary) representation of the group $G$ on the Hilbert
space $\Cal H$ is a map $\pi : G \rightarrow \Cal U(\Cal H)$
satisfying $\pi(g)\pi(h)= \mu_{g,h} \pi(gh), \forall g,h \in G$, for
some scalar 2-{\it cocycle} $\mu$ on $G$, i.e. $\mu : G \times G
\rightarrow \Bbb T$ satisfies $\mu_{g,h}
\mu_{gh,k}=\mu_{h,k}\mu_{g,hk}$, $\forall g,h,k \in G$. It is
immediate to see that the equivalence class of such a $\pi$ only
depends on the class of $\mu$ in H$^2(G,\Bbb T)\overset \text{\rm
def} \to = \text{\rm Z}^1(G,\Bbb T)/\text{\rm B}^1(G,\Bbb T)$, where
Z$^1(G,\Bbb T)$ denotes the multiplicative group of all scalar
valued 2-cocycles and B$^1(G,\Bbb T)$ is the subgroup of {\it
coboundaries}, $\mu_{g,h}=\lambda_g \lambda_h
\overline{\lambda_{gh}}$, for some $\lambda:G \rightarrow \Bbb T$.

\proclaim{1.1. Lemma} Let  $H \subset G$ be an 
inclusion of groups satisfying the relative property
$\text{\rm (T)}$. Fix $1\geq \varepsilon > 0$ and let  
$F(\varepsilon), \delta(\varepsilon)$ 
be the constants given by $(1.0)$. Denote $\tilde
F(\varepsilon)=F(\varepsilon^2/28)$, $\tilde
\delta(\varepsilon)=\delta(\varepsilon^2/28)/2$. 
Then the following holds true:

\vskip .05in

If $\pi :
G\rightarrow \Cal{U}(\Cal{H})$ is a projective representation with
scalar 2-cocycle $\mu\in \text{\rm H}^2(G,\Bbb T)$, and $\xi \in
\Cal{H}$ is a unit vector satisfying $d(\pi
(g)\xi,\Bbb{C}\xi)\leq\tilde \delta(\varepsilon), \forall g\in
\tilde F(\varepsilon)$, then $\exists \xi_0\in\Cal{H}$ and
$\lambda:H\rightarrow \Bbb{T}$, such that
$\Vert\xi-\xi_0\Vert<\varepsilon$, $\pi(h)\xi_0=\lambda_h\xi_0$ and
$\mu_{h,h'}=\lambda_h\lambda_{h'}\overline{\lambda}_{hh'}$, $\forall
h,h'\in H$.

In particular, if $\delta_1=\frac{1}{2}\delta(\frac{1}{28})$,
$F_1=F(\frac{1}{28})$, then whenever  $\pi:G\rightarrow\Cal{U(H)}$
is a projective representation with scalar 2-cocycle $\mu$ such that
$\|\pi(g)\xi-\xi\|<\delta_1,\forall g\in F_1$ for some unit vector
$\xi\in\Cal{H}$, $\mu_{|H}$ follows coboundary.
\endproclaim
\noindent {\it Proof}.
 Note first that if $\pi: G\rightarrow\Cal U(\Cal H)$ is a projective
representation of the group $G$ on the Hilbert space $\Cal H$, then
$\pi\otimes\overline\pi: G\rightarrow\Cal U(\Cal H\otimes \overline
{\Cal H} )$ is a genuine representation of $G$ on the Hilbert space
$\Cal H\otimes \overline {\Cal H}$.
Identify $\Cal H\otimes \overline {\Cal H}$ with $\Cal{HS}$, the
Hilbert space of Hilbert-Schmidt operators  on $\Cal H$,  and then
note that $\pi\otimes\overline\pi$ can be extended to all $\Cal
B(\Cal H)$ by the formula
$(\pi\otimes\overline\pi)(g)(T)=\pi(g)T\pi(g)^*$ for all $g\in G$
and for all $T$ in $\Cal B(\Cal H)$.

Fix $\varepsilon>0$ and let  $\tilde F(\varepsilon)$ and $\tilde
\delta(\varepsilon)$ be defined as in the second part of the
statement. Let $\xi\in \Cal H$ be a unit vector such that $d(\pi
(g)\xi,\Bbb{C}\xi)\leq\tilde \delta(\varepsilon), \forall g\in
\tilde F(\varepsilon)$. Then for all $g\in
F(\varepsilon^2/28)$:$$\|(\pi\otimes\overline\pi)(g)(\xi\otimes\overline\xi)-
\xi\otimes\overline\xi\|=\|\pi(g)\xi\otimes\overline\pi(g)\overline\xi-
\xi\otimes\overline\xi\|_\Cal{HS}<\delta(\varepsilon^2/28)$$

By the relative property $\text{\rm(T)}$ applied to the
representation $\pi\otimes\overline\pi$ on $\Cal {HS}$, there exist
a Hilbert-Schmidt operator $T$ of $\Cal{HS}$-norm equal to $1$ such
that $(\pi\otimes\overline\pi)(h)(T)=T,\forall h\in H$ and
$\|T-\xi\otimes\overline\xi\|\leq\|T-
\xi\otimes\overline\xi\|_\Cal{HS}<\varepsilon^2/28.$


Thus, $\pi(h)T\pi(h)^*=T, \forall h\in H$, implying that the
operators $\pi(h)$ and $T$ on the Hilbert space $\Cal H$ commute,
$\forall h\in H$. But then  $T^*$ and $TT^*$ also commute with
$\pi(h), \forall h\in H$. Thus, all the spectral projections of
$TT^*$ are in the commutant of $\pi(H)$ in $\Cal B(\Cal H)$.  Since
$TT^*$ is a trace class operator, its spectral projections have
finite trace, i.e. they are finite dimensional.

Since $\|TT^*-\xi\otimes\overline\xi\|<2\varepsilon^2/28$, it
follows that $\|(TT^*)^2-\xi\otimes\overline\xi\|<4\varepsilon^2/28$
and  $\|(TT^*)^2-TT^*\|<6\varepsilon^2/28$. Thus, if
$\varepsilon<1$, then  there exist a non-zero spectral projection
$P$ of $TT^*$ with finite rank such that $\|P-TT^*\|<
12\varepsilon^2/28$. This implies that
$\|P-\xi\otimes\overline\xi\|<\varepsilon^2/2$.

In particular $P$ has to be a rank one projection, i.e. $P$ is of
the form $\xi'\otimes\overline\xi'$ and $|\langle
\xi',\xi\rangle|>1-\varepsilon^2/2$. Taking $\alpha \in \Bbb T$ such
that $\alpha \langle \xi', \xi \rangle >0$ we get that
$\|\alpha\xi'-\xi\|<\varepsilon$.

 Let $\xi_0=\alpha\xi'$. Then $\pi(h)\xi_0 \otimes
 \overline{\pi(h)(\xi_0)}=
 (\pi\otimes\overline\pi)(h)(\xi_0\otimes\overline\xi_0)
 =\xi_0\otimes\overline\xi_0$
for all $h\in H$ which implies that $\pi(h)(\xi_0) =\lambda_h\xi_0$
with $\lambda_h\in \Bbb T$. Since $\pi (h)\pi(h')\xi_0=\mu_{h,h'}
\pi(hh')\xi_0$ we get $\lambda_h\lambda_{h'}=
\mu_{h,h'}\lambda_{hh'}$ for all $h,h'\in H$. \hfill $\square$

\proclaim{1.2. Theorem} Let $M$ be a separable finite von Neumann
algebra. Let $G$ be a discrete group with a subgroup $H$ such that
$(G,H)$ has the relative property $\text{\rm (T)}$. Let
$\{\pi_j\}_{j \in J}$ be projective representations of $G$ into the
unitary group of $p_jMp_j$, with scalar $2$-cocycles
$\{\mu_{_j}\}_{j\in J}$, where $p_j\in \Cal P(M)$ are non-zero
projections in $M$. Then the image of the set $\{\mu_j|_{H}\}_{j\in
J}$ in $\text{\rm H}^2(H,\Bbb T)$ is at most countable.

\endproclaim
\noindent {\it Proof}.  Let $J_0 \subset J$ be such that the
cocycles $\mu_j|_H, j\in J_0$, are distinct in $\text{\rm
H}^2(H,\Bbb T)$ and $\{\mu_j|_H\}_{j\in J_0}=\{\mu_j|_H\}_{j\in J}$.
We have to prove that $J_0$ is countable.

Assume it is not. Then there exists $c>0$ such that the set $J_1=\{j
\in J_0 | \tau(p_j) \geq c\}$ is uncountable. Let $F_1$ and
$\delta_1$ be as in Lemma 1.1. Let also $\tau$ be a normal
faithful trace state on $M$ and denote as usual by $\|x\|_2 =
\tau(x^*x)^{1/2}$ for $x\in M$ the corresponding Hilbert norm on
$M$. Since $M$ is separable and $J_1$ is uncountable, there exist
$j_1,j_2\in J_1$ such that:

$$
\|p_{j_1}-p_{j_2}\|_2 < \delta_1c/4
$$
$$
\Vert \pi_{j_1}(g)-\pi_{j_2}(g) \Vert_{2} < \delta_1c/4, \forall
g\in F_1.
$$

In particular, the first inequality shows that

$$
\|p_{j_1}\|_2 - \|p_{j_1}p_{j_2}\|_2 \leq
\|p_{j_1}(p_{j_1}-p_{j_2})\|_2 \leq \delta_1c/4 \leq c/2,
$$
implying that $\|p_{j_1}p_{j_2}\|_2 \geq c - c/2 = c/2$.

For $x\in M$, denote by $L(x),R(x)$ the operators of left,
respectively right multiplication by $x$ on $L^2(M,\tau)$. Define
$\pi:G\longrightarrow \Cal{B}(p_{j_1}L^2(M,\tau)p_{j_2})$ by
$\pi(g)\eta=L(\pi_{j_1}(g))R(\pi_{j_2}(g)^*)\eta$. Then $\pi$ is a
projective representation of cocycle $\mu_{j_1}\bar{\mu}_{j_2}$ and
if we denote
$\xi=\|p_{j_1}p_{j_2}\|_{2}^{-1}(p_{j_1}p_{j_2})^{\hat{}}$ then
$\xi$ has norm one and we have for all $g \in F_1$ the estimates:

$$\align \|\pi(g)\xi-\xi\| & =\|p_{j_1}p_{j_2}\|_{2}^{-1}\|\pi_{j_1}(g)p_{j_2}-
p_{j_1}\pi_{j_2}(g)\|_2\\& = \|p_{j_1}p_{j_2}\|_{2}^{-1}
\|(\pi_{j_1}(g)-\pi_{j_2}(g))p_{j_2}+(p_{j_2}-
p_{j_1})\pi_{j_2}(g)\|_2
\\& \leq(c\delta_1/4+c\delta_1/4)\|p_{j_1}p_{j_2}\|_{2}^{-1}\leq (c\delta_1/2)
\|p_{j_1}p_{j_2}\|_{2}^{-1} \leq \delta_1  \endalign
$$

>From lemma 1.1 it follows that the cocycle
$\mu_{j_1}\bar{\mu}_{j_2}$ is a coboundary in $H$, which contradicts
$\mu_{j_1}|_{H}\not=\mu_{j_2}|_{H}$. \hfill $\square$ \vskip .1in

For the next corollary, recall that given any scalar 2-cocycle $\mu
$ on a discrete group $G$, one associates to it the {\it projective
left regular representation} $\lambda_\mu : G \rightarrow \Cal
U(L^2(G))$, with scalar 2-cocycle $\mu$, defined by $\lambda_\mu(g)
(\Sigma_h c_h \xi_h)=\Sigma_h c_h \mu_{g,h} \xi_{gh}$. We denote
$L_\mu(G)=\lambda_\mu(G)''$ the corresponding von Neumann algebra,
as considered by Connes and Jones in ([CJ]).

\proclaim{1.3. Corollary} Let $H\subset G$ be an inclusion of
discrete groups with the relative property $\text{\rm (T)}$. Let $M$
be a separable finite von Neumann algebra. Let $J$ be the set of
scalar $2$-cocycles $\mu \in \text{\rm H}^2(G,\Bbb T)$ such that
$L_\mu(G)$ can be embedded into $M$ $($not necessarily unitaly$)$.
Then the set $\{\mu_{|H} \mid \mu \in J\}\subset \text{\rm
H}^2(H,\Bbb T)$ is countable.
\endproclaim
\noindent {\it Proof}.  It is enough to show that for every $n$ the
set: $\{\mu_{|H}\mid \mu \in \text{\rm H}^2(G,\Bbb T)$ with
$L_\mu(G)$ embeddable (not necessarily unitaly) in $M_n(M)\}$ is
at most countable. Since for any scalar 2-cocycle $\mu$ for $G$ the
$\mu$-twisted left regular representation $\lambda_\mu$ is a
projective representation with scalar 2-cocycle $\mu$ and whose von
Neumann algebra is $L_\mu(G)$, the statement follows from the
previous theorem. \hfill $\square$ \vskip .1in

Note that when applied to the case $M=M_{n\times n} (\Bbb C)$ the
proof of lemma 1.1 gives an estimate of the number of certain sets
of scalar 2-cocycles of $H$ in terms of the constants of rigidity of
$H \subset G$. We emphasize this in the next proposition, where we
also include an estimate of the number of projective representations
of dimension $n$ of a group with the property $\text{\rm(T)}$,
generalizing a result in ([HRV]).

We need the following notations: For $G$ a discrete group and $n$ a
positive integer, we denote by $\Cal{PR}(G,n)$ the set of
equivalence classes of projective representations of $G$ of
dimension $n$. Also, we denote by $\text{\rm H}^2(G,n)$ the set of
scalar 2-cocycles $\mu\in\text{\rm H}^2(G,\Bbb T)$ for which there
exists a projective representation $\pi\in\Cal{PR}(G,n)$ with
cocycle $\mu$.

\proclaim{1.4. Proposition} $1^\circ$. Let $G$ be a discrete group with
property $\text{\rm (T)}$. There exists a constant $c>1$, which
depends only on the constants of rigidity of G, such that
$\Cal{PR}(G,n)$ has at most $c^{n^2}$ elements, $\forall n\geq 1$.

$2^\circ$. Let $H\subset G$ be an inclusion of discrete groups such
that $(G,H)$ has the relative property $\text{\rm(T)}$. There exists
a constant $d$ such that the subset $\{\mu_{|H} \mid \mu\in\text{\rm
H}^2(G,n)\}$ of $\text{\rm H}^2(H,\Bbb T)$ has at most $d^{n^2}$
elements, $\forall n\geq 1$.

\endproclaim
\noindent {\it Proof}.  $1^\circ$. Let $(F_1,\delta_1)$ be as in
Lemma 1.1, for $H=G$. By [Wy] there exists $c_0>1$ such that one can
cover the unit sphere in $\Bbb{R}^{2n^2}$ with $c_0^{n^2}$ balls of
radius $\delta_1 /2$, $\forall n\geq 1$. Thus one can cover the
sphere of radius $\sqrt{n}$ with $c_0^{n^2}$ balls of radius
$\sqrt{n}\delta_1 /2$.

We first show that there are at most $c_0^{|F_1|n^2}$ irreducible
projective representations of dimension $n$. Assume not. Regarding
unitaries in $M_n(\Bbb{C})$ as vectors of norm $\sqrt{n}$ in $\Bbb
R^{2n^2}$, the pigeonhole principle implies that there exist
$\pi_1,\pi_2$ irreducible such that: $\Vert \pi_1(g)-\pi_2(g)
\Vert_2<\delta_1,\forall g\in F_1$.

As in the proof of Theorem 1.2, define $\pi:G\rightarrow
\Cal{B}(L^2(\Bbb{M}_n(\Bbb{C}),\tau))$ by
$\pi(g)\eta=L(\pi_{1}(g))R(\pi_{2}(g)^*)\eta$, $\forall \eta \in
L^2(M_{n\times n}(\Bbb{C}),\tau)$, where $\tau$ is the normalized trace.
Then $\pi$ is a projective representation of $G$ with cocycle
$\mu_\pi = \mu_{\pi_1} \overline{\mu}_{\pi_2}$ and it satisfies
$\Vert \pi(g)\hat{1}-\hat{1} \Vert_2 <\delta_1, \forall g\in F_1$.
Lemma 1.1 implies that there exists $\lambda:G\rightarrow \Bbb{T}$
and a unit vector $\xi_0\in\Bbb{M}_n(\Bbb{C})$ such that if we
define $\pi'_1(g)= \overline\lambda_g\pi_1(g)$ then
$\pi'_1(g)\xi_0=\xi_0\pi_2(g), \forall g,$ and $\pi_1', \pi_2$ have
the same cocycle $\mu'$. Taking the adjoint and noticing that
$\pi'_1(g)\pi'_1(g^{-1})=\pi_2(g)\pi_2(g^{-1})=\mu'_{g, g^{-1}}1,
\forall g\in G$, it follows $\xi_0^*\pi'_1(g)=\pi_2(g)\xi_0^*,
\forall g\in G$. The two relations imply that $\xi_0^*\xi_0$
commutes with $\pi_2(g), \forall g\in G$, and using the
irreducibility of $\pi_2$ it follows that $\pi'_1,\pi_2$ are
conjugate.

Thus, for every $n$ there are at most $c_0^{|F_1|n^2}$ irreducible
projective representations of order $n$. This implies that the
number of projective representations of dimension $n$ is at most
$2^nc_0^{|F_1|n^2}$, so $c=2c_0^{|F_1|}$ will do.

$2^\circ$. Let $(F_1,\delta_1)$ be the constants from Lemma 1.1 for
$(G,H)$, and let $c_0$ be as before. Let $d=c_0^{|F_1|}$. Assume
that $\{\mu_{|H}, \mid \mu\in\text{\rm H}^2(G,n)\}$ has more than
$d^{n^2}$ elements. By the pigeonhole principle if follows that
there exist $\pi_{j_1},\pi_{j_2}\in\Cal{PR}(G,n)$ with cocycles
$\mu_{j_1},\mu_{j_2}$ such that $\mu_{j_1}|_{H}\not=\mu_{j_2}|_{H}$
and $\Vert \pi_{j_1}(g)-\pi_{j_2}(g) \Vert_2<\delta_1, \forall g\in
F_1$. This leads to a contradiction, as in the proof of Lemma 1.1.
\hfill $\square$

\vskip .1in \noindent {\bf 1.5. Remark.} The same proof as for part
$1^\circ$ of the above proposition shows that the similar result for
groups $G$ with the property $(\tau)$ of Lubotzky holds true (see
1.3 in [LZ] for the definition of property $(\tau)$ and 1.4.3 in
[LZ] for related statements).

\heading {2. Examples} \endheading

Recall that the groups $SL(n, \Bbb Z)$ $n>2$ and $Sp(2n,\Bbb Z)$
$n>1$ have the property $\text{\rm(T)}$ of Kazhdan, \cite{Ka}.
Obvious examples of inclusions with relative property
$\text{\rm(T)}$ are $H\subset H\times\Gamma$ with $H$ a property
$\text{\rm(T)}$ group and $\Gamma$ an arbitrary discrete group.
It is shown in ([Ka], [Ma])
that $\Bbb
Z^2\subset \Bbb Z^2\rtimes SL(2,\Bbb Z)$ has the relative property
$\text{\rm(T)}$. More generally, by a result of Burger \cite{Bu},
any inclusion of the form
$\Bbb Z^2\subset \Bbb Z^2\rtimes \Gamma$, with $\Gamma
\subset SL(2,\Bbb Z)$ non-amenable, has the relative property
$\text{\rm(T)}$. Recently Valette \cite{Va} showed that if $\Gamma $
is an arithmetic lattice in an absolutely simple Lie group, then
there exist an embedding of $\Gamma $ in $SL(m,\Bbb Z)$ for some
$m$, such that $\Bbb Z^m\subset \Bbb Z^m\rtimes \Gamma$ has the
relative property $\text{\rm(T)}$. Fernos \cite{Fe} constructed
other examples of
inclusions of groups $\Bbb Z^m \subset \Bbb Z^m \rtimes \Gamma$
with the relative property $\text{\rm(T)}$, for $\Gamma\subset GL(m,\Bbb Z)$.

More examples of pairs of groups having the relative property (T)
come out from the following easy observation:

\proclaim{2.1. Lemma} Let $\sigma: \Gamma \rightarrow \text{\rm
Aut}(H)$ and $\sigma': \Gamma \rightarrow \text{\rm Aut}(H')$ be
actions of a $\Gamma$ on $H, H'$ and denote by $\tilde{\sigma}:
\Gamma \rightarrow \text{\rm Aut}(H \times H')$ the diagonal action,
$\tilde{\sigma}(g)(x,y)=(\sigma(g)x, \sigma'(g)y)$, $x\in H, y\in
H'$, $g \in \Gamma$.

$1^\circ$. If $H \subset H \rtimes_\sigma \Gamma$ and $H' \subset
H'\rtimes_{\sigma'} \Gamma$ have the relative property $\text{\rm
(T)}$ then $(H \times H') \subset (H \times H')
\rtimes_{\tilde{\sigma}} \Gamma$ has the relative property
$(\text{\rm T})$.

$2^\circ$. Assume $H \subset H \rtimes \Gamma$ has the relative
property $\text{\rm (T)}$. Let $\beta\in \text{\rm Aut}(\Gamma)$.
Denote $\sigma'=\sigma \circ \beta$ and $\tilde{\sigma}$ the
diagonal action $\tilde{\sigma}(g)=\sigma(g) \times \sigma'(g)$ of
$G$ on $H \times H$. Then $(H \times H) \subset (H \times H)
\rtimes_{\tilde{\sigma}} \Gamma$ has the relative property
$(\text{\rm T})$.

$3^\circ$. If  $\Gamma $ is a subgroup of $GL(n,\Bbb Z)$
 such that $\Bbb Z^n\subset \Bbb Z^n\rtimes \Gamma$ has the
relative property $\text{\rm(T)}$ then $\Bbb Z^{2n}\subset \Bbb
Z^{2n}\rtimes_\theta \Gamma$ also has the relative property
$(\text{\rm T})$, where for each $g\in\Gamma$ $\theta(g)\in
SL(2n,\Bbb Z)$ is the matrix $\pmatrix
g&0\\
0&(g^{-1})^t
\endpmatrix $.

\endproclaim
\noindent {\it Proof.} $1^\circ$. If $\pi$ is a unitary
representation of $(H \times H') \rtimes \Gamma$ on $\Cal H$ with an
almost invariant unit vector $\xi\in \Cal H$ then by $(1.0)$ $\xi$
follows uniformly almost invariant to $H \times \{e'\}$ and to
$\{e\} \times H'$. Since
$$
\|\pi(h,h')\xi - \xi\|\leq \|\pi(h,e')\pi(e,h') \xi - \pi(h,e') \xi
\|+\|\pi(h,e')\xi-\xi\| \
$$
$$
\leq \|\pi(e,h') \xi -  \xi \| + \|\pi(h,e')\xi-\xi\|
$$
for all $h\in H$, $h'\in H'$, it follows that $\xi$ is uniformly
almost invariant to $H \times H'$. Thus, if $\xi_0$ is the element
of minimal norm in $\overline{\text{\rm co}}^{w} \{ \pi(h,h')\xi;
h\in H,h'\in H'\} \subset \Cal H$, then $\xi_0$ is invariant to
$\pi(H \times H')$ and $\xi_0 \neq 0$.

$2^\circ$. Since the inclusions $H \subset H \rtimes_{\sigma}
\Gamma$ and $H \subset H \rtimes_{\sigma'} \Gamma$ are isomorphic,
and the first inclusion has the relative property $\text{\rm(T)}$,
the second one has this property  as well. Thus, part 1$^\circ$
applies to get the conclusion.

$3^\circ$. Apply $2^\circ$ to $(H \subset H \rtimes \Gamma$)=$(\Bbb
Z^n\subset \Bbb Z^n \rtimes\Gamma)$ and $\beta(g)=(g^{-1})^t$.
 \hfill $\square$
\vskip .2in

Of all these examples of inclusions of groups $H \subset G$ with the
relative property (T) we are interested in those for which the set
of restrictions of 2-cocycles $\{\mu|_H:\,\,\mu\in \text{\rm
H}^2(G,\Bbb T)\}$ is ``large'' (uncountable), so we can take
advantage of Corollary 1.3.
  There are difficulties in obtaining such examples.
First it is difficult to calculate second cohomology groups.
Secondly it is hard to control the size of this group when
restricted to $H$. We overcome these difficulties by looking at
inclusions of the form $(H\subset G)=(\Bbb Z^{2n}\subset \Bbb
Z^{2n}\rtimes \Gamma)$, and the 2-cocycles on $G$ arise as
 extensions to $G$ of $\Gamma$ invariant 2-cocycles in $H$.
A similar construction has been considered in \cite{Cho2}.

\vskip .1in
 Denote with $J$ the matrix
 $\pmatrix
 0 & I_n \cr
-I_n & 0
  \endpmatrix\in GL(2n,\Bbb Z)$.
It defines a 2-cocycle $\nu:\Bbb Z^{2n}\times\Bbb Z^{2n}\to \Bbb Z$
by the formula $\nu(x,y)=x^tJy$. For each $\alpha\in \Bbb T$ we
denote with $\nu_\alpha$ the $\Bbb T$-valued 2-cocycle  defined by
 $\nu_\alpha\overset \text{\rm def}\to =\alpha^{\frac{1}{2}\nu}$.
 Since $\nu_\alpha$
 is a coboundary iff $\alpha=1$, $\alpha\mapsto\nu_\alpha$
 is an embedding of $\Bbb T$ into
 $\text{\rm H}^2(\Bbb Z^{2n},\Bbb T)$.

The set of invertible matrices that leave $\nu$ (and also
$\nu_\alpha$) invariant is the symplectic group  $Sp(2n,\Bbb Z)$.
Thus, given any subgroup $\Gamma$ of $Sp(2n,\Bbb Z)$, $\nu_\alpha$
can be extended to a  2-cocycle on $\Bbb Z^{2n}\rtimes\Gamma$, which
we still denote $\nu_\alpha$,
 by the formula
$\nu_\alpha((x_1,\gamma_1),(x_2,\gamma_2))\overset \text{\rm def}\to
=\nu_\alpha (x_1,\gamma_1 x_2).$

\vskip .1in \noindent {\bf 2.2. Notations}. For each subgroup
$\Gamma \subset Sp(2n,\Bbb Z)$ and each $\alpha\in \Bbb T$ let
$N_\alpha$ and $M_\alpha(\Gamma,n)$ be the cocycle von Neumann
algebras
 $N_\alpha\overset \text{\rm def}\to =
L_{\nu_\alpha}(\Bbb Z^{2n})\subset L_{\nu_\alpha}(\Bbb
Z^{2n}\rtimes\Gamma) \overset \text{\rm def} \to =
M_\alpha(\Gamma,n)$. Alternatively, $M_\alpha(\Gamma,n)$ can be
regarded as the cross product von Neumann algebra
$N_\alpha\rtimes_{\sigma_\alpha}\Gamma$, where the action
$\sigma_\alpha$ is defined by
 $\sigma_\alpha(g)(\lambda_{\nu_\alpha}(x))=\lambda_{\nu_\alpha}(gx)$ for all
 $x\in \Bbb Z^{2n}$ and all $g\in \Gamma$. Note that the isomorphism
class of $M_\alpha(\Gamma,n)$ may in fact depend on the embedding
$\Gamma \subset Sp(2n,\Bbb Z)$, in other words it may depend on the
way $\Gamma$ acts on $\Bbb Z^{2n}$, a fact that is not well
emphasized by the notation $M_\alpha(\Gamma,n)$. For instance, the
group $\Bbb F_2$ can be embedded in $SL(2,\Bbb Z)$ in many ways,
giving different actions of $\Bbb F_2$ on $\Bbb Z^2$ and thus on
$L_{\nu_\alpha}(\Bbb Z^{2})$.

\vskip .1in

If $\alpha$ is a root of unity of order $m$ then $N_\alpha$ is
homogeneous of type $\text{\rm I}_{nm}$, while if $\alpha$ is not a
root of unity, then $N_\alpha$ is isomorphic to the hyperfinite
$\text{\rm II}_1$ factor $R$. \noindent $N_\alpha' \cap M_\alpha
(\Gamma,n) =\Cal Z(N_\alpha)$. Also, if either $\Bbb Z^{2n}$ has no
$\Gamma$-invariant finite subsets other than $\{0\}$ or if $\alpha$
is not a root of unity then $M_\alpha(\Gamma,n)$ is a $\text{\rm
II}_1$ factor.

In the case when $n=1$, $Sp(2,\Bbb Z)$ is in fact equal to
$SL(2,\Bbb Z)$ and if we denote by $u$ and $v$ the canonical
generators of $\Bbb Z^2$ we have that
$\lambda_u\lambda_v=\nu_\alpha(u,v)\lambda_{uv}=
\alpha^\frac{1}{2}\lambda_{uv}$ and
$\lambda_v\lambda_u=\nu_\alpha(v,u)\lambda_{vu}=
\alpha^{-\frac{1}{2}}\lambda_{uv}$ showing that
$\lambda_u\lambda_v=\alpha\lambda_v\lambda_u$. So when $n=1$ and
$\alpha=\exp(2\pi\imath \theta)$ with $\theta$ irrational,
$L_{\nu_\alpha}(\Bbb Z^2)$ is the hyperfinite $\text{\rm II}_1$
factor represented as the irrational rotation algebra $R_\alpha$ of
angle $\theta$. Thus if $\Gamma$ is an arbitrary non amenable
subgroup of $SL(2,\Bbb Z)$, then $M_\alpha(\Gamma)\overset \text{\rm
def}\to = M_\alpha(\Gamma,1)$ is an irrational rotation HT factor,
as considered in (\cite{P2},\cite{P4}).

 Recall that two finite von Neumann algebras $M$ and $N$
are {\it stably isomorphic} if $M$ is isomorphic to an
amplification $N^t$ of $N$, i.e., if there exist $n\in \Bbb{N}$ and
a projection $p\in M_n(N)$ such that $M$ is isomorphic to $p
M_n(N)p$ ($= N^t$, where $t = n\tau(p)$).

 \proclaim{2.3. Corollary}
 Let $M$ be a separable $\text{\rm II}_1$ factor. If $\Gamma$ is
 a subgroup of $Sp(2n,\Bbb Z)$ such that $\Bbb
Z^{2n}\subset\Bbb Z^{2n}\rtimes\Gamma$ has the relative property
$(\text{\rm T})$,
 then the
set of $\alpha \in \Bbb T$ for which some amplification of
$M_\alpha(\Gamma,n)$ can be embedded into $M$ is at most countable.
Thus, the factors $\{ M_{\alpha}(\Gamma,n)\}_{\alpha \in \Bbb T}$
 are non stably isomorphic modulo countable sets.
\endproclaim

\proclaim{2.4. Corollary}
$1^\circ.$
 The irrational rotation HT factors
$M_\alpha(\Gamma)$ cannot be all embedded into a separable
$\text{\rm II}_1$ factor and are non stably isomorphic modulo
countable sets.

$2^\circ.$
If  $\Gamma $ is a subgroup of $GL(n,\Bbb Z)$
 such that $\Bbb Z^n\subset \Bbb Z^n\rtimes \Gamma$ has the
relative property $\text{\rm(T)}$ then the factors
$M_\alpha(\Gamma,n)=L_{\nu_\alpha}( \Bbb Z^{2n}\rtimes_\theta
\Gamma)$ where $\theta$ is as in $2.1.3^\circ$ cannot be all
embedded into a separable $\text{\rm II}_1$ factor and are non
stably isomorphic modulo countable sets.

\endproclaim

\vskip .15in

\head 3. Disintegration of rigid von Neumann algebras \endhead

We now use results from the previous section to derive some
properties of the disintegration of the type $\text{\rm II}_1$ von
Neumann algebras coming from property $\text{\rm(T)}$ groups
$\Lambda$ with large (infinite) radical, i.e. for which $L(\Lambda)$
has diffuse center. Thus, we give an explicit description of the
disintegration of the type II$_1$ von Neumann algebras from the
property (T) groups $\Lambda$ with large center constructed by Serre
(see \cite{dHVa} pg. 40), which arise as central extensions of
property (T) groups. We also use an argument from ([O1]) to show
that the factors in the disintegration of the algebra $L(\Lambda)$
of an arbitrary property (T) group $\Lambda$ are mutually
non-isomorphic, modulo countable sets.

We first recall some facts about the disintegration theory of a von
Neumann algebra (see chapter 2 of \cite {D1} for a detailed
treatment). Thus, let $(\Cal Z,\mu)$ be a Borel space with a
positive measure and $\xi\to\Cal H_\xi$ be a measurable field of
Hilbert spaces on $\Cal Z$. We denote by $\Cal H=\int_{\xi\in\Cal
Z}^\oplus \Cal H_\xi \text{ }d\mu$ the
 corresponding direct integral Hilbert space.
An operator field $\xi\to T_\xi$,  $T_\xi\in\Cal B(\Cal H_\xi)$ is
called diagonalizable if it is of the form $\xi\to c(\xi)I_{\Cal
H_\xi}$ where $c:\Cal Z\to \Bbb C$ is measurable. An operator $T$
acting on  $\Cal H$ is called decomposable if it comes from a
measurable operator field $\xi\to T_\xi$, in which case we write
$T=\int_{\xi\in\Cal Z}^\oplus T_\xi \text{ }d\mu$. An operator $T$
is decomposable if and only if it commutes
with the set of diagonalizable operators.

Now assume that for each $\xi\in\Cal Z$, $\Cal A_\xi$ is a von
Neumann algebra acting on $\Cal H_\xi$. $\xi\to\Cal A_\xi$ is a
measurable field of von Neumann algebras if there exist a sequence
$\{T_i\}_{i\in \Bbb Z}$ of measurable operator fields such that for
each
 $\xi\in\Cal Z$, $\{T_i(\xi)\}_{i\in\Bbb Z}$
generates $\Cal A_\xi$. The set of decomposable operators
$T=\int_{\xi\in\Cal Z}^\oplus T_\xi \text{ }d\mu$ for which
$T_\xi\in\Cal A_\xi$ is a von Neumann algebra and it is denoted by
 $\Cal A=\int_{\xi\in\Cal Z}^\oplus \Cal A_\xi \text{ }d\mu$.

\vskip .1in \noindent {\it 3.0 Example:} Let $G$ be a discrete group
with a 2-cocycle $\nu:G\times G\rightarrow A$ where $A$ is a
discrete abelian group. The central extension of $G$ with cocycle
$\nu$ is a group $\tilde{G}$ where $\tilde{G}=A\times G$ as a set
and the multiplication is given by $(a_1,g_1)(a_2,g_2)=(a_1a_2\nu
(g_1,g_2),g_1g_2)$. Notice that $(a_1,g_1)^{-1}=(a_1^{-1}\nu
(g_1,g_1^{-1})^{-1},g_1^{-1})$. By a result of Serre, if $G$ is a
property $\text{\rm(T)}$ group and $\nu\not=0$ in $\text{\rm
H}^2(G,A)$ then $\tilde{G}$ also has property $\text{\rm(T)}$.

 For each character $\alpha\in \hat A$
let $L_\alpha (G)=L_{\nu_\alpha}(G)$, where
 $\nu_\alpha$ is the $\Bbb T$-valued 2-cocycle given by the formula
$ \nu_\alpha(g_1,g_2)=\alpha(\nu(g_1,g_2))$.

 Let $B=C^*_{red}(\tilde{G})$ and $\tau$ be the natural
trace on $B$ defined by $\tau(a,g)=\delta_{(a,g)}^{(e_A,e_G)}$.
   For each $\alpha\in\hat A$ let $\tau_\alpha$ be the trace on $B$
   defined by $\tau_\alpha(a,g)=\alpha(a)\delta_g^{e_G}$.
    Let $(\pi, \Cal{H_\circ})$ and $(\pi_\alpha,\Cal{H}_\alpha)$ be the GNS
    representations of $B$ with respect to the states $\tau$ and $\tau_\alpha$.
     Then: $\Cal{H}_\circ=\ell^2(\tilde{G})$
     , $\Cal{H}_\alpha=\ell^2(G)$, $\pi(B)''=L(\tilde{G})$ and
     $\pi_\alpha (B)''\simeq L_\alpha(G)$.
      The last equality is easy to check since
      $\tau_\alpha((a-\alpha(a))(a-\alpha(a))^*)=0$ so $\pi_\alpha(a)=\alpha(a) I$.

For each $g\in G$ define the vector field $x_g$ to be $x_g(\alpha)=
\widehat{(e_A,g)}^{\Cal H_\alpha}$, where  for any $b\in B$,
$\widehat{\,b\,}^{\Cal H_\alpha}$ denotes the class of $b$ in $\Cal
H_\alpha$.
 It is clear that for each
$\alpha\in \hat A$ fixed, the set $\{x_g(\alpha)\}_{g\in G}$ is an
orthonormal basis of $\Cal H_\alpha$ and that for each $g_1,g_2\in
G$ the function $\alpha\to \langle
x_{g_1}(\alpha),x_{g_2}(\alpha)\rangle_{\Cal H_\alpha}$
 is continuous.
Then by \cite{D1} II,1,4 proposition 4, there exist a unique
structure of measurable Hilbert spaces on
 $\alpha\mapsto\Cal H_\alpha$ that make the vector fields $x_g$ measurable.
  Moreover, a vector field $x$ is measurable
 if and only if
$\alpha\to \langle x_g(\alpha),x(\alpha)\rangle_{\Cal H_\alpha}$ is
measurable for every $g\in G$.

 Let $\theta:\Cal H_\circ \rightarrow \int_{\alpha\in\hat
A}^\oplus \Cal H_\alpha \text{ }d\alpha$
 be the linear map defined by:
 $$\theta\widehat{(a,g)}^{\Cal H_\circ} =(\widehat{(a,g)}^{\Cal H_\alpha})_{\alpha\in\hat A}
 =(\alpha(a) g)_{\alpha\in\hat A}$$
\noindent
 (where the last equality is via the identification $\Cal{H}_\alpha=l^2(G)$.)

  We show that $\theta$ is an isomorphism of Hilbert spaces and:
  $$(\theta(\pi(x)\xi))_\alpha=\pi_\alpha(x)\theta(\xi), \forall x\in B, \xi\in \Cal H_\circ$$

 Note that

\noindent
 $$\langle\theta\widehat{(a_1,g_1)}^{\Cal H_\circ},
\theta\widehat{(a_2,g_2)}^{\Cal H_\circ}\rangle_{\int^\oplus \Cal H_\alpha }=
 \int_{\alpha\in\hat A}\langle\widehat{(a_1,g_1)}^{\Cal H_\alpha},\widehat{(a_2,g_2)}^{\Cal H_\alpha}\rangle
 _{\Cal H_\alpha}\text{ }d\alpha=$$
 $$
 =\int_{\alpha\in\hat A}\tau_\alpha(a_1a_2^{-1}\nu(g_1,g_2^{-1})\nu(g_2,g_2^{-1})^{-1},g_1g_2^{-1})\text{ }d\alpha,$$
 with the last term being zero whenever $(a_1,g_1)\neq
(a_2,g_2)$. Thus
$$\langle\widehat{\theta(a_1,g_1)},\widehat{\theta(a_2,g_2)}\rangle_{\int^\oplus \Cal H_\alpha d\alpha}=
\delta_{(a_1,g_1)}^{(a_2,g_2)}=\langle\widehat{(a_1,g_1)}^{\Cal
H_\circ},\widehat{(a_2,g_2)}^{\Cal H_\circ}\rangle_{\Cal H_\circ}$$
showing that
 $\theta$ is an injective morphism of
Hilbert spaces.

 To check surjectivity, let $\{x(\alpha)\}_{\alpha}\in
\int^\oplus \Cal H_\alpha \text{ }d\alpha$ be
 a measurable vector field.
  Since for every fixed $g\in G$ the function $\alpha\mapsto \langle x(\alpha), x_g(\alpha)\rangle$
   belongs to $L^2(\hat A)$, there exist $(d_{a,g})_{a\in A, g\in G}$
    such that: $\sum_{a\in A} d_{a,g}\alpha(a)= \langle x(\alpha), x_g(\alpha)\rangle
    , \forall \alpha\in \hat A, g\in G$. Moreover
     $\sum_{a,g} |d_{a,g}|^2=\int \Vert x(\alpha) \Vert ^2 \text{ }d\alpha$
     is finite. Define $v=\sum_{a,g} d_{a,g}(a,g)$. Then $v\in \Cal H_\circ$ and
     $\theta (v)=\{x(\alpha)\}_\alpha$, which shows that $\theta$ is an  isomorphism.
\vskip .05in
 We now check that
$(\theta(\pi(x)\xi))_\alpha=\pi_\alpha(x)\theta(\xi)_\alpha, \forall x\in
B,\text { }\xi\in \Cal H_\circ$. This is clear since
$(\theta(\pi(x)\xi))_\alpha=\widehat{(\pi(x)\xi)}^{\Cal H_\alpha}=
\widehat{(x\xi)}^{\Cal H_\alpha}=\pi_\alpha(x)\widehat{(\xi)}^{\Cal
H_\alpha}=\pi_\alpha(x)\theta(\xi)_\alpha.$

The diagonalizable operator fields on $\int_{\alpha\in\hat A}^\oplus \Cal H_\alpha d\alpha$ correspond, via $\theta^{-1}$, to the elements of the von Neumann algebra $\pi (A)''\subset B(\Cal H_0)$. Altogether, using 8.4.1 in \cite{D2} we have thus obtained:

\proclaim{3.1. Proposition} Let $G$ be a discrete group with a
2-cocycle $\nu:G\times G\rightarrow A$, where $A$ is a discrete
abelian group. Let $\tilde{G}$ be the central extension of $G$
 defined by the cocycle $\nu$. For each $\alpha\in \hat{A}$
 let $L_\alpha (G)=L_{\nu_\alpha}(G)$, where
 $\nu_\alpha$ is the $\Bbb T$-valued 2-cocycle on $G$ defined as
$ \nu_\alpha(g_1,g_2)=\alpha(\nu(g_1,g_2))$. Then the von Neumann
algebra $L(\tilde{G})$
 has the following direct integral decomposition:

$$L(\tilde{G})=\int_{\alpha\in\hat A}^\oplus L_\alpha (G)\text{ } d\alpha$$
\endproclaim
\vskip .2in

 Note that if we let $G=\Bbb Z^{2n}\rtimes \Gamma$, $n>2$,
 where $\Gamma$ is a non amenable subgroup of $Sp(2n,\Bbb Z)$
and $\nu$ the $\Bbb Z$-valued 2-cocycle on $G$ defined in Section 2,
then from Corollary 2.2 and Proposition 3.1 it follows that the
factors in the above direct integral decomposition of $L(\tilde G)$
are property $\text{\rm(T)}$ and they are non-isomorphic modulo
countable sets.

But in fact one can obtain a general result along these lines, by
using an argument similar 
to Ozawa's proof that there are 
no ``universal'' 
separable II$_1$ factors ([O1]). We include the details of 
the argument, 
for completeness.

\proclaim{3.2. Theorem} Let $\Lambda$ be a discrete property
$\text{\rm(T)}$ group such that the von Neumann algebra $L(\Lambda)$
has diffuse center, and let $L(\Lambda)=\int_{t\in\Cal Z}^\oplus M_t
\text{ }d\mu$ be its direct integral decomposition. Then there
exists a set $\Cal Z_0\subset \Cal Z$, $\mu(\Cal Z_0)=0$, such that
the factors $M_t$, $t\in\Cal Z\setminus\Cal Z_0$, are mutually
non-stably isomorphic modulo countable sets.
\endproclaim

\noindent {\it Proof.} Let $B=C^*_{red}(\Lambda)$,
let $\tau$ be the canonical trace of $B$ and let $\Cal Z=\widehat{Z(B)}$.
 The direct integral decomposition of the GNS representation of $(B,\tau)$
  induces factorial representations $\pi_t:\Lambda\rightarrow \Cal B (\Cal H_t)$, $t\in \Cal Z$.
   The factors in the direct integral decomposition of $L(\Lambda)$ are
    $M_t=\pi_t(\Lambda)''\subset \Cal B(\Cal H_t)$,
    and we may assume $\Cal H_t= L^2(M_t)$.
    By 8.4.1 and 8.4.2 in \cite{D2}
    there exists a measure zero set $\Cal Z_0\subset \Cal Z$ such that the representations
     $\pi_t$, $t\in \Cal Z\setminus\Cal Z_0$, are mutually non-conjugate.

    Assume, by contradiction, that $M_t$ is isomorphic to an amplification $M^{s(t)}$
  of the same factor $M$, for all $t\in S$, where $S\subset \Cal Z\setminus\Cal Z_0$ is uncountable.
   We may clearly assume $c \leq s(t)\leq 1$, $\forall t$, for some
   $c>0$. To simplify notations,
   we still denote by $\pi_t$ the representations of $\Lambda$
   into the unitary group of $p_tMp_t$,
   induced by the isomorphisms $M_t\simeq p_tMp_t$,
   where $p_t\in \Cal P(M)$, $\tau(p_t)=s(t)$, $t\in S$.

Let $(F_0,\delta_0)$ be property (T) constants for $\Lambda$ as defined in section 1.
 By using a separability argument
 as in Theorem 1.2 and ([GP]),
 it follows that there exist $t_1\not=t_2\in S$ such that $p_{t_1}$
 is close to $p_{t_2}$ and such that
     if  $\pi:\Lambda\rightarrow \Cal B(p_{t_1}L^2(M)p_{t_2})$
     denotes
     the representation  of $\Lambda$ given
     by the formula
   $\pi(g)\eta=L(\pi_{t_1}(g))R(\pi_{t_2}(g)^*)\eta$
    and $\xi$ is the vector $\xi=\Vert p_{t_1}p_{t_2}\Vert ^{-1}(p_{t_1}p_{t_2})^{\hat{}}$
    then
    $\Vert \pi(g) \xi - \xi\Vert_2<\delta_0$ for all $g\in F_0$.
Since $\Lambda$ has property (T), there exists a non-zero vector
$\eta\in p_{t_1}L^2(M)p_{t_2}$ such that $\pi(g)\eta=\eta$, for all
$g\in \Lambda$.
   Equivalently, if we regard $\eta$ as a square integrable operator,
   we have $\pi_{t_1}(g)\eta=\eta\pi_{t_2}(g)$, for all $g\in \Lambda$.
   By the standard trick, if $v\in M$ is the partial isometry in the polar decomposition of $\eta$
   with the property that the right supports of $\eta$ and $v$ coincide,
   then $vv^*\in \pi_{t_1}(\Lambda)'\cap p_{t_1}Mp_{t_1}=\Bbb Cp_{t_1}$,
   $v^*v\in \pi_{t_2}(\Lambda)'\cap p_{t_2}Mp_{t_2}=\Bbb Cp_{t_2}$ and
   $\pi_{t_1}(g)v=v\pi_{t_2}(g)$, for all $g\in \Lambda$.
   This implies that
   $\pi_{t_1},\pi_{t_2}$ are conjugate representations of $\Lambda$,
   which contradicts $t_1\not=t_2$.


 \
\hfill $\square$

\vskip .5in

\heading 4. Conjugacy and isomorphism problems for
$M_{\alpha}(\Gamma)$ \endheading

We have seen that the cocycle von Neumann algebras
$M_{\alpha}(\Gamma)$ constructed in Section 2 can be regarded as the
crossed product von Neumann algebras
$R_\alpha\rtimes_{\sigma_\alpha}\Gamma$. Moreover, by ([P4]), when
$\alpha\in \Bbb T$ is irrational the isomorphism class of the
algebras $M_\alpha(\Gamma)$ is completely determined by the cocycle
conjugacy class of the actions $\sigma_\alpha$ of $\Gamma$ on the
hyperfinite II$_1$ factor $R \simeq R_\alpha$. Thus, the
classification of the factors $M_\alpha(\Gamma)$ amounts to the
classification up to cocycle conjugacy of the actions
$(\sigma_\alpha, \Gamma)$. In particular, for a fixed $\Gamma\subset
SL(2,\Bbb Z)$, showing that the factors $M_\alpha(\Gamma)$ are
non-isomorphic for different irrational numbers $\alpha$ amounts to
showing that the corresponding actions $\sigma_\alpha$ are
non-cocycle conjugate. While we cannot solve this latter problem, we
show here that for a large class of subgroups $\Gamma\subset SL(2,
\Bbb Z)$ the conjugacy class of the action $\sigma_\alpha$
determines the irrational number $\alpha$.

\proclaim{4.1. Theorem} Let $\Gamma \subset SL(2,\Bbb{Z})$ be a
subgroup of $SL(2,\Bbb{Z})$ containing a parabolic element $a$
 and an element $b$ that doesn't commute with $a$. If
$\alpha_1$ and $\alpha_2$ are irrationals in the upper half torus
such that the actions $\sigma_{\alpha_1}$ and $\sigma_{\alpha_2}$ of
$\Gamma$ on the hyperfinite $\text{\rm II}_1$ factors
$R_{\alpha_j}=L_{\mu_{\alpha_j}}(\Bbb Z^2)$ ($j=1,2$)
 are conjugate then $\alpha_1=\alpha_2$.
\endproclaim
\noindent {\it Proof}. By replacing $\Gamma$ with
$\gamma\Gamma\gamma^{-1}$ for a certain $\gamma\in SL(2,\Bbb{Z})$,
we may assume that $a$ has $(1,0)$ as eigenvector. We may also
assume that the corresponding eigenvalue is 1, by substituting $a$
with $a^2$ if necessary. Thus $a=\pmatrix 1 & n \cr 0 & 1 \cr \endpmatrix \in \Gamma$, for some $n\in \Bbb{Z}$ non-zero.

For $j=1,2$ let $\alpha_j=e^{2\pi it_j}, \alpha_j^\frac{1}{2}=e^{\pi
it_j}$ with $t_j\in [0,1/2) \setminus\Bbb Q$, and let
$u_j=\lambda_{\mu_{\alpha_j}}(1,0)$, and
$v_j=\lambda_{\mu_{\alpha_j}}(0,1)$ be the unitaries generating
$L_{\mu_{\alpha_j}}(\Bbb{Z}^2)=R_{\alpha_j}$. The cocycle relation
$u_jv_j=\alpha_j v_ju_j$ implies $u_j^kv_j^l=\alpha_j^{kl}
v_j^lu_j^k$ for all $k,l\in \Bbb{Z}$.
 For $g=\pmatrix p & q \cr r & s \cr
\endpmatrix \in \Gamma$ we have:
%
%
$$\sigma_{\alpha_j}(g)(u^kv^l)=\alpha_j^{\frac{1}{2}(kl-(pk+ql)(rk+sl))}u_j^{pk+kl}v_j^{rk+sl}$$

Assume $\sigma_{\alpha_1}$ and $\sigma_{\alpha_2}$ are conjugate,
i.e. there exists an isomorphism $\theta:R_{\alpha_{2}}\rightarrow
R_{\alpha_{1}}$ such that
$\theta(\sigma_{\alpha_2}(g)(x))=\sigma_{\alpha_1}(g)(\theta(x))$ for all $g \in \Gamma, x\in R_{\alpha_{2}} $. We prove
$\alpha_1=\alpha_2$.

 Denote $u=u_1,v=v_1,u'=\theta(u_2),v'=\theta(v_2)$. To
simplify notations, we identify $x\in R_{\alpha_{1}}$ with its image
$\hat{x}$ in $L^2(R_{\alpha_{1}})$. Thus
$(u^kv^l)_{(k,l)\in\Bbb{Z}^2}$ is an orthonormal basis of
$L^2(R_{\alpha_1},\tau)$ and $R_{\alpha_1}$ is identified with the
set of ``Fourier expansions'' $\sum_{(k,l)\in
\Bbb{Z}^2}\lambda_{k,l} u^kv^l$ in $L^2(R_{\alpha_1},\tau)$, that are
(twisted) left convolvers on $L^2(R_{\alpha_1},\tau)$. Let
$$u'=\sum_{(k,l)\in \Bbb{Z}^2}c_{k,l}
u^kv^l\;,\;\;\;\; v'=\sum_{(k,l)\in \Bbb{Z}^2}d_{k,l} u^kv^l$$ for
some $c_{k,l},d_{k,l}\in \Bbb{C}$ such that $$\sum_{(k,l)\in \Bbb{Z}^2}
\vert {c_{k,l}} \vert ^2 < \infty\;,\;\;\;\; \sum_{(k,l)\in
\Bbb{Z}^2} \vert {d_{k,l}} \vert ^2 < \infty$$
%
%
Since the actions $\alpha_1,\alpha_2$ are conjugate via $\theta$, we have:
$$
\sigma_{\alpha_1}(g)((u')^k(v')^l)=
\alpha_2^{\frac{1}{2}(kl-(pk+ql)(rk+sl))}(u')^{pk+ql}(v')^{rk+sl}
$$
\noindent
Choosing  $g=a, k=1, l=0$, we obtain  $\sigma_{\alpha_1}(a)(u')=u'$. Thus

$$
\sigma_{\alpha_1}(a)(\sum_{(k,l)\in \Bbb{Z}^2}c_{k,l}u^kv^l)
=\sum_{(k,l)\in \Bbb{Z}^2}c_{k,l}\alpha_1^{-\frac{1}{2}nl^2}
u^{k+nl}v^l=\sum_{(k,l)\in \Bbb{Z}^2}c_{k,l}u^kv^l
$$
which implies
$\alpha_1^{-\frac{1}{2}nl^2}c_{k-nl,l}=c_{k,l}, \forall k,l\in
\Bbb{Z}$. Thus, for $l$ nonzero: $\vert c_{k,l}\vert =\vert
c_{k-nl,l}\vert =\vert c_{k-2nl,l} \vert =...$ have to be all zero
since $\sum_{(k,l)\in \Bbb{Z}^2} \vert {c_{k,l}} \vert ^2 <
\infty$. Denote $c_k=c_{k,0}$. Then $u'=\sum_{k\in
\Bbb{Z}}c_{k}u^k$.

Let $b=\pmatrix m_1 & m_2 \cr m_3 & m_4 \cr \endpmatrix,
m_1m_4-m_2m_3=1$. $ab\not =ba$ is equivalent to $m_3\not =0$.
Using the formula for $\sigma_{\alpha_1}$ for $g=b,k=1,l=0$,
we obtain $\sigma_{\alpha_1}(b)(u')=\alpha_2^{-\frac{1}{2}m_1m_3}
(u')^{m_1}(v')^{m_3}$.

This implies $u'\sigma_{\alpha_1}(b)(u') =\alpha_2^{-\frac{1}{2}m_1m_3}
(u')^{m_1}u'(v')^{m_3} ={\alpha_2}^{m_3}\sigma_{\alpha_1}(b)(u')u'$,

\vskip .1in
\noindent
 and thus
$$
\sum_{k\in \Bbb{Z}}c_{k} u^k(\sum_{j\in \Bbb{Z}}c_{j}
\alpha_1^{-\frac{1}{2}m_1m_3j^2}u^{m_1j}v^{m_3j})
$$
$$
={\alpha_2}^{m_3}(\sum_{j\in \Bbb{Z}}c_{j}
\alpha_1^{-\frac{1}{2}m_1m_3j^2}u^{m_1j} v^{m_3j})\sum_{k\in
\Bbb{Z}}c_{k} u^k
$$
Hence we obtain:

$$
\sum_{k,j\in \Bbb{Z}}c_{k}c_{j}\alpha_1^{-\frac{1}{2}m_1m_3j^2}
(1-{\alpha_2}^{m_3}{\alpha_1}^{-m_3kj}) u^{k+m_1j}v^{m_3j}=0.
$$

\noindent Since the function $(k,j)\rightarrow (k+m_1j,m_3j)$ is injective for
$m_3\not =0$, it follows:

$$c_{k}c_{j}({\alpha_1}^{m_3kj}-{\alpha_2}^{m_3})=0,
\forall k,j\in\Bbb{Z}$$

Letting  $k=j$ we obtain $c_k=0$, for all $k$ except possibly
two values $k_0,-k_0$. Indeed, since $\alpha_1$ is not a root of
unity there exists at most one $N=m_3k^2$ such that
${\alpha_1}^{N}={\alpha_2}^{m_3}$.

Since $u'$ is not a scalar, we know $k_0\not=0$. Taking $j=-k_0$ and
using
${\alpha_1}^{-m_3k_0^2}\not={\alpha_1}^{m_3k_0^2}={\alpha_2}^{m_3}$
we obtain $c_{k_0}c_{-k_0}=0$. Thus only one coefficient of the
Fourier expansion of $u'$ is non-zero.  So far we have showed then
that:
$$u'=cu^{k_0}\text{ and } {\alpha_1}^{k_0^2}={\alpha_2}.$$
Now substituting $u'$ in the relation $u'v'=\alpha_2 v'u'$ we
obtain:

$$
cu^{k_0}(\sum_{(k,l)\in \Bbb Z^2}d_{k,l}u^kv^l)
=\alpha_2(\sum_{(k,l)\in \Bbb Z^2}d_{k,l}u^kv^l)cu^{k_0}
$$
thus

$$
\sum_{(k,l)\in \Bbb Z^2}d_{k,l}u^{k_0+k}v^l =\sum_{(k,l)\in \Bbb
Z^2}d_{k,l}\alpha_2\alpha_1^{-k_0l}u^{k+k_0}v^l
$$
which yields

$$
d_{k,l}(1-\alpha_2\alpha_1^{-k_0l})=0, \forall k,l\in \Bbb{Z}
$$

Since $\alpha_1^{k_0l}\not=\alpha_2$ unless $l=k_0$ we obtain that
$d_{k,l}=0$, for all $k\in\Bbb{Z}$ and $l\not=k_0$. Denote
$d_k=d_{k,k_0}$. Thus we have $ u'=cu^{k_0}$ and $v'=(\sum_{k}
d_ku^k)v^{k_0}. $ This implies that for every $j\geq 1$
there exists $w_j\in W^{*}(1,u)$ such that $(v')^j=w_jv^{jk_0}$.
 Using the formula for $\sigma_{\alpha_1}$
one more time for $g=b$, $k\not=0$ arbitrary and $l=1$, we have:
\vskip.1in

 $ \sigma_{\alpha_1}(b)((u')^kv')
=\alpha_2^{\frac{1}{2}(k-(m_1k+m_2)(m_3k+m_4))}
(u')^{m_1k+m_2}(v')^{m_3k+m_4}=$
 \vskip .1in
$=\alpha_2^{\frac{1}{2}(k-(m_1k+m_2)(m_3k+m_4))}c^{m_1k+m_2}
u^{k_0(m_1k+m_2)}w_{m_3k+m_4}v^{k_0(m_3k+m_4)}$

\vskip .1in

 \noindent On the other hand:  $\sigma_{\alpha_1}(b)((u')^kv')
=\sigma_{\alpha_1}(b)(\sum_{l}c^{k}d_lu^{kk_0+l}v^{k_0})= $
\vskip.1in
\noindent
 $ \sum_{l}c^{k}d_l\alpha_1^{\frac{1}{2}[(kk_0+l)k_0-(m_1(kk_0+l)+m_2k_0)
(m_3(kk_0+l)+m_4k_0)]} u^{m_1(kk_0+l)+m_2k_0}v^{m_3(kk_0+l)+m_4k_0}
$ \vskip .1in Identifying the corresponding coefficients, for every
$l$ we must have either $d_l=0$ or
$m_3(kk_0+l)+m_4k_0=k_0(m_3k+m_4)$, which implies $l=0$. Thus
$d_l=0,\forall l\not=0$ and $v'=dv^{k_0}$ for some $d\in\Bbb{C}$.
Altogether, $u'=cu^{k_0}, v'=dv^{k_0}$ for some $c,d\in \Bbb C$.
Since $u',v'$ generate $R_{\alpha_1}$, this implies $k_0=1$
or $k_0=-1$. But $\alpha_1^{-1}\not=\alpha_2$ because
$\alpha_1,\alpha_2$ belong to the upper half torus.
Thus $k_0=1$  and $\alpha_1=\alpha_2$. \hfill $\square$ \vskip .1in

\heading APPENDIX: A general result on fundamental groups
\endheading

We give here a short proof of a result in ([P2]), showing that the
HT factors $M_\alpha(\Gamma)$, $\alpha \in \Bbb T$, $\Gamma \subset
SL(2,\Bbb Z)$ non-amenable, have at most countable fundamental
group. The result we prove is in fact much more general, covering
all results of this type in ([P1,2]), as particular cases:

\proclaim {A.1. Theorem} Let $M$ be a separable $\text{\rm II}_1$
factor. Assume there exists a non-zero projection $p\in M$ such that
$pMp$ contains a von Neumann subalgebra $B$ such that $B \subset
pMp$ is a rigid inclusion and $B'\cap pMp \subset B$. Then $\Cal
F(M)$ is countable.
\endproclaim
\noindent {\it Proof}. Recall from (4.2 of [P2]) that $B \subset M$
rigid implies there exist $F \subset M$ finite and $\delta > 0$ such
that if $\phi:M \rightarrow M$ is a subunital, subtracial completely
positive map which satisfies $\|\phi(x)-x\|_2 \leq \delta, \forall
x\in F$, then $\|\phi(u)-u\|_2 \leq 1/2, \forall u\in \Cal U(B)$.

Since the fundamental groups of $M$ and $pMp$ coincide, it is
clearly sufficient to prove the statement in the case $p=1$. For
each $t \in (0, 1) \cap \Cal F(M)$ choose a projection $p_t \in \Cal
P(B)$  and an isomorphism

$\theta_t : M \simeq p_tMp_t$. Since $B$ is diffuse we can make the
choice so that in addition we have $p_t \leq p_{t'}$ whenever  $t
\leq t'$.

Assume $\Cal F(M)$ is uncountable. Thus, $[c,1) \cap \Cal F(M)$ is
uncountable for some $0< c < 1$. By the separability of $M$, this
implies there exist $ t, s \in \Cal F(M) \cap [c,1)$, $t < s$,  such
that $\|\theta_s(x)-\theta_t(x)\|_2 \leq \delta c$, $\forall x\in
F$.

Thus, if we denote $\theta = \theta_s^{-1} \circ \theta_t$ then

$\theta$ is an isomorphism of $M$ onto $qMq$, where $q =
\theta_s^{-1} (\theta_t (1))= \theta_s^{-1}(p_t)$, and we have
$\theta(1) \leq 1$, $\tau(q) \geq c$, $\tau \circ \theta \leq \tau$,
$\|\theta(x)-x\|_2 \leq  \delta$, $\forall x\in F$. Consequently, we
have $\|\theta(u)-u\|_2 \leq 1/2, \forall u \in \Cal U(B)$.

Let $k$ denote the unique element of minimal norm $\|\quad\|_2$ in
$K=\overline{\text{\rm co}}^w \{ \theta(u)u^*\mid u\in \Cal U(B)\}$.
Then $\|k-1\|_2\leq 1/2$ and thus $k\neq 0$. Also, since
$\theta(u)Ku^* \subset K$ and $\|\theta(u)ku^*\|_2 = \|k\|_2,
\forall u \in \Cal U(B)$, by the uniqueness of $k$ it follows that
$\theta(u)ku^*=k$, or equivalently $\theta(u)k = ku$, for all $u\in
\Cal U(B)$. By a standard trick, if $v\in M$ is the (non-zero)
partial isometry in the polar decomposition of $k$ and if we express
any element in $B$ as linear combination of unitaries, then we get
$\theta(b)v=vb, \forall b\in B$, $v^*v \in B'\cap M = \Cal Z(B),
vv^* \in \theta(B)'\cap qMq = \Cal Z(\theta(B)q)$.

Since in particular $v^*v\in B$, we can apply the above to $b=v^*v$
to get $\theta(v^*v)v=v v^*v$. But this implies
$\theta(v^*v)vv^*=vv^*$, so that $\theta(v^*v)\geq vv^*$. This is a
contradiction, since $\theta$ shrinks the trace of any elements by
$\tau (q)<1$, while $\tau(vv^*)=\tau(v^*v)$. \hfill $\square$

\vskip .1in \vskip .05in

\proclaim {A.2. Corollary} For each $\Gamma \subset SL(2,\Bbb Z)$
non-amenable and $\alpha \in \Bbb T$, the factor $M_\alpha(\Gamma)$,
as defined in Sections $0$ and $2$, has countable fundamental group.
\endproclaim
\noindent {\it Proof}. Since $\Bbb Z^2 \subset \Bbb Z^2 \rtimes
\Gamma$ has the relative property $\text{\rm(T)}$ (cf [Bu]), the
inclusion of von Neumann algebras $R_\alpha=L_\alpha(\Bbb Z^2)
\subset L_\alpha(\Bbb Z^2 \rtimes \Gamma)=M_\alpha(\Gamma)$ has the
relative property $\text{\rm(T)}$ and $R_\alpha'\cap
M_\alpha(\Gamma)\subset R_\alpha$. Thus A.1 applies. \hfill
$\square$

\head  References\endhead

\item{[B]} B. Brenken: {\it Representations and
automorphisms of the irrational rotation algebra}, Pacific Journal
of Mathematics, {\bf 111} (1984), 257-282.

\item{[Bu]} M. Burger:
{\it Kazhdan constants for} $SL(3,\Bbb Z)$, J. Reine Angew. Math.,
{\bf 413} (1991), 36-67.

\item{[Ch1]} M. Choda: {\it A continuum of non-conjugate
property} $\text{\rm(T)}$ {\it actions of} $SL(n,\Bbb{Z})$ {\it on
the hyperfinite} $\text{\rm II}_1$ {\it factor}, Math. Japon., {\bf
30} (1985), 133-150.

\item{[Ch2]} M. Choda: {\it Outer actions of
groups with property} $\text{\rm(T)}${\it on the hyperfinite}
$\text{\rm II}_1$ {\it factor }, Math. Japon., {\bf 31} (1986),
533-551.

\item{[C]} A. Connes: {\it A factor of type $\text{\rm II}_1$ with countable
fundamental group}, J. Operator Theory, {\bf 4} (1980) p. 151-153.

\item{[CJ]} A. Connes, V.F.R. Jones: {\it Property} (T) {\it for
von Neumann algebras}, Bull. London Math. Soc. {\bf 17} (1985),
57-62.

\item{[D1]} J. Dixmier: ``Les alg\`ebres
d'op\'erateurs dans l'espace Hilbertien'',
Gauthier-Villars 1957.

\item{[D2]} J. Dixmier: ``Les C* alg\`ebres et leurs repr\'esentations'',
Gauthier-Villars 1964.

\item {[Fe]} T. Fernos: {\it Kazhdan's relative property $\text{\rm
(T)}$: Some New Examples}, math.GR/0411527

\item{[GP]} D. Gaboriau and S. Popa: {\it An uncountable
family of non orbit equivalent actions of $\Bbb{F}_n$}, J. Am. Math.
Soc. {\bf 18} (2005), 547-559; math.GR/0306011.

\item{[H]} U. Haagerup: {\it An example of non-nuclear C$^*$-algebra
which has the metric approximation property}, Invent. Math. {\bf 50}
(1979), 279-293.

\item{[dHVa]} P. de la Harpe, A. Valette:  La propri\'et\'e $\text{\rm(T)}$ de
Kazhdan pour les groupes localement compacts, Ast\'erisque {\bf 175}
(1989)

\item{[HRVa]} P. de la Harpe, A. G. Robertson,
A. Valette: {\it On the spectrum of the sum of generators for a
finiteley generated group}, Israel J. Math. {\bf 81} (1993), 65-96.

\item{[Jo]} P. Jolissaint: {\it On Property} (T)
{\it for pairs of topological groups}, l'Ens. Math. {\bf 51} (2005),
31-45.

\item{[LZ]} A. Lubotzky, A. Zuk: {\it On the property
$(\tau)$}, preprint 2002.

\item{[M]} G. Margulis: {\it Finitely-additive invariant measures
on Euclidean spaces}, Ergodic Theory Dynam. Systems, {\bf 2} (1982),
383-396.

\item{[N]} O. A. Nielsen: ``Direct integral theory'', Lecture
Notes in Pure and Appl. Math. {\bf 61}, Marcel Dekker, New York,
1980.

\item{[O1]} N. Ozawa: {\it There is no separable universal
$\text{\rm II}_1$-factor}, Proc. Amer. Math. Soc., {\bf 132} (2004), 487-490.
math.OA/0210411

\item{[O2]} N. Ozawa:
{\it A Kurosh type theorem for factors}, math.OA/0401121.

\item{[P1]} S. Popa: {\it Correspondences}, INCREST preprint 56/1986,
unpublished.

\item{[P2]} S. Popa: {\it On a class of type} $\text{\rm II}_1$
{\it factors with Betti numbers invariants}, Ann. of Math., {\bf
163} (2006), 809-899. math.OA/0209310.

\item{[P3]} S. Popa: {\it Some rigidity results for
non-commutative Bernoulli shifts}, J. Fnal. Analysis {\bf 230}
(2006), 273-328.

\item{[P4]} S. Popa: {\it A unique
decomposition result for HT factors with torsion free core},
\newline math.OA/0401138.

\item{[R]} M. Rieffel: {\it $C^*$-algebras
associated with irrational rotations},
Pacific. J. Math. {\bf 93} (1981), 415-429.

\item{[Sh]} Y. Shalom: {\it Bounded generation and Kazhdan's property} $\text{\rm(T)}$,
Publ. Math. I.H.E.S. {\bf 90}(2001), 145-168.

\item {[Va]} A. Valette: {\it Group pairs with relative property
$\text{\rm (T)}$ from arithmetic lattices}, preprint 2004
(preliminary version 2001).

\item{[V]} D. Voiculescu: {\it Property} $\text{\rm(T)}$ {\it
and approximation of operators},
Bull. London Math. Soc {\bf 22} (1990), 25-30.

\item{[Wy]} A. D. Wyner:
{\it Random packings and coverings of the unit n-sphere}, Bell
System Tech. J., {\bf 46} (1967), 2111-2118.

\enddocument